**The differential role of verbal and visuospatial working memory in mathematics and reading**

David Giofrè[1], Enrica Donolato[2], & Irene C. Mammarella[2]

[1] *Natural Science and Psychology, Liverpool John Moores University, Liverpool (UK)*

[2] *Department of Developmental and Social Psychology, Padova (Italy)*

Correspondence concerning this article should be addressed to:

Dr. David Giofrè
Liverpool John Moores University
Natural Sciences and Psychology
Tom Reilly Building Byrom Street, Liverpool, L3 3AF
Tel. +44 151 904 6336
Fax. +44 151 904 6302
E-mail: david.giofre@gmail.com

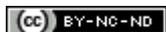

## Abstract

**Objectives**: Several studies have focused on the role of working memory (WM) in predicting mathematical and reading literacy. Alternative models of WM have been proposed and a modality-dependent model of WM, distinguishing between verbal and visuospatial WM modalities, has been advanced. In addition, the relationship between verbal and visuospatial WM and academic achievement has not been extensively and consistently studied, especially when it comes to distinguishing between mathematical and reading tasks. **Method:** In the present study, we tested a large group of middle school children in several measures of WM, and in mathematical and reading tasks. **Results:** Confirmatory factor analyses showed that verbal and visuospatial WM can be differentiated and that these factors have a different predictive power in explaining unique portions of variance in reading and mathematics. **Conclusions:** Our findings point to the importance of distinguishing between WM modalities in evaluating the relationship between mathematics and reading.

## Research highlights

- WM tasks can be distinguished between verbal and visuospatial modalities
- Verbal WM predicts a large portion of unique variance in reading
- Visuospatial WM predicts a large portion of unique variance in mathematics

**The differential role of verbal and visuospatial working memory on mathematics and reading**

Academic success is an essential aspect of everyday life for children and is also particularly important to future occupational outcomes (Schmidt & Hunter, 2004). The literature suggests that several different cognitive abilities are involved in academic achievement. In particular, working memory (WM) has long been considered one of the most important predictors of academic performance (St Clair-Thompson & Gathercole, 2006). Variations in children's WM capacities are related to progress across many areas of academic achievement including mathematics and reading literacy (Gathercole, Pickering, Knight, & Stegmann, 2004; Jarvis & Gathercole, 2003). In the present study, the role of WM in verbal and visuospatial components, was taken into account to better understand WM's involvement in mathematics and reading tasks.

**1.1 The influence of WM on mathematical and reading tasks**

WM is a limited-capacity system that enables information to be temporarily stored and manipulated (Baddeley, 2000), and it has repeatedly been associated with academic performance (Cornoldi & Giofrè, 2014). Alternative WM models have been proposed. The most classical WM model, the tripartite model (Baddeley & Hitch, 1974), includes a central executive, responsible for controlling resources and monitoring information, and two domain-specific modalities for either verbal or visuospatial information. Other researchers have suggested a modality-independent model, according to which WM is seen as a domain-general factor, and the distinction is made between a WM factor, which requires cognitive control to a large extent, and a short-term memory (STM) factor, which requires less cognitive control (i.e., less attentional resources) (Kane et al., 2004). A domain-specific factors model, only distinguishing between verbal and visuospatial modalities, has also been proposed (Shah & Miyake, 1996). Finally, others maintain the view that WM is best represented by a single model (Pascual-Leone, 1970).

Several independent studies with children demonstrate that the classical tripartite model is superior compared to alternative WM models (Gathercole, Pickering, Ambridge, & Wearing, 2004; Giofrè, Borella, & Mammarella, 2017; Giofrè, Mammarella, & Cornoldi, 2013). Intriguingly, there is a debate on whether or not is possible to distinguish between modalities (verbal and visuospatial) in complex working memory tasks (i.e., tasks that require high levels of cognitive control). On the one hand, some authors claimed that this distinction is not possible because the correlation between verbal (WM-V) and visuospatial (WM-VS) working memory tasks is too strong (Kane et al., 2004). Conversely, other authors supported a distinction between WM-V and WM-VS in certain age groups (Swanson, 2008). Nevertheless, understanding the structure of WM is important because of its relationship with educational outcomes and academic achievement (Cornoldi & Giofrè, 2014).

There is extensive evidence of the importance of WM in key academic domains such as reading (Gathercole, Alloway, Willis, & Adams, 2006), reading comprehension (Borella & de Ribaupierre, 2014; Oakhill, Hart, & Samols, 2005), school readiness (Fitzpatrick & Pagani, 2012), and mathematical achievement (Berg, 2008; Bull, Espy, & Wiebe, 2008; De Smedt et al., 2009; Passolunghi, Mammarella, & Altoè, 2008; Swanson & Beebe-Frankenberger, 2004), including mental addition (Caviola, Mammarella, Cornoldi, & Lucangeli, 2012; Mammarella et al., 2013), mental subtraction (Caviola, Mammarella, Lucangeli, & Cornoldi, 2014), geometry (Giofrè, Mammarella, & Cornoldi, 2014; Giofrè, Mammarella, Ronconi, & Cornoldi, 2013), and problem solving (Passolunghi & Mammarella, 2010; Swanson & Sachse-Lee, 2001).

However, most of the previously cited research focused on individual difference measures, and only few studies have considered the relationship between verbal and visuospatial WM in mathematics or reading achievement mesures for children with typical development. For example, several studies show that verbal WM is related to reading (Swanson, Xinhua Zheng, & Jerman, 2009), and visuospatial WM is related to math skills (Caviola et al., 2014; Holmes, Adams, & Hamilton, 2008), whereas other research suggested the same WM resources can be attributed to both reading and mathematics (e.g., Swanson, 1999, 2011). Additionally, Wilson and Swanson (2001) found that both verbal WM and visuospatial WM predicted mathematics achievement, but that verbal WM accounted for most of the significant variance in predictions of math performance.

In a recent study (Swanson, 2017), in which participants between the ages of 5 and 90 were tested, it was shown that after controlling for age, WM accunted for a large portion of the variance in both reading and mathematics, respectively. Also, verbal and visuospatial tasks seem to predict some unique variance in reading (7% by WM-VS and 6% by WM-V) and in mathematics (4% by WM-VS and 2% by WM-V) in children (Swanson, 2017). This finding seems to indicate that WM-VS and WM-V have as similar predictive power on mathematics and reading. However, several pieces of evidence seem to indicate that reading is more influenced by verbal WM (Daneman & Carpenter, 1980; Swanson et al., 2009), while mathematics is more influenced by visuospatial WM (Caviola et al., 2014; Mammarella, Caviola, Giofrè, & Szűcs, 2017; Szűcs, Devine, Soltesz, Nobes, & Gabriel, 2013).

In the present study, we aimed to test whether or not verbal and visuospatial WM processes can be separated in middle school children with typical development in middle school children. The relationship between WM and academic achievement has been extensively studied in primary school children (Alloway & Alloway, 2010; Alloway et al., 2005; Bull et al., 2008; Gathercole & Pickering, 2000), while less attention has been devoted to students attending middle schools (e.g., Gathercole et al., 2004). Moreover, if this can be done, the study aimed to examine whether or not the WM modality is important in predicting performance in reading (which is supposedly more influenced by verbal WM) and in mathematics (believed to be more influenced by visuospatial WM).

## 2. Method

### 2.1 Participants

The study involved 144 middle schoolchildren in sixth and eighth grades. Three participants were found to be multivariate outliers using Cook's distance and were excluded from the analyses, so the final sample included a total of 141 children (71 males and 70 females, $M_{age}$ = 149.89 months, SD = 13.54) in sixth grade (N = 66, females = 53%) and eighth grade (N = 75, females = 47%). The children came from middle-class families (87.0% Italian, 7.0% African; 2.1% Chinese; 3.9% Other), and were attending schools in an urban area of the north-east Italy.

The study was approved by the Ethics Committee on Psychology Research at the University of Padova (Italy). Parental consent was obtained. Children with special educational needs, intellectual disabilities, or neurological/genetic disorders were not included in the study.

**2.2 Materials**

**2.2.1 Working memory tasks**

**Verbal WM**

*Word span -Backward* (WB-B; Cornoldi & Vecchi, 2003). This task can be considered as a complex working memory task in children (Alloway, Gathercole, & Pickering, 2006; Gathercole et al., 2006). Words were presented verbally at a rate of 1 item per second, proceeding from the shortest series to the longest (from 2 to 8 items). There was no time limit for recalling the words in the same forward order (Cronbach's alphas .75). The score was the number of words accurately recalled in the correct order.

*Dual tasks verbal (DT-V)* (De Beni, Palladino, Pazzaglia, & Cornoldi, 1998). The DT-V material consisted of orally-presented word lists, each containing four words of high-medium frequency. The word lists were grouped into sets of various sizes (i.e., containing from 2 to 6 word lists each). The children were asked to press the space bar whenever they heard an animal noun and, after completing each set, they had to recall the last word on each list, in the same order of presentation. This task showed adequate psychometric properties (Cronbach's α = .69) and a good predictive power (Giofrè, Mammarella, & Cornoldi, 2013).

*Listening span test* (LST; Daneman & Carpenter, 1980; Palladino, 2005). The children listened to sentences grouped into sets of various sizes (containing from 2 to 5 sentences each).

After hearing each sentence, children were asked whether the sentence was true or false. After completing each set, the children had to recall the last word in each sentence, in the order of presentation. The score corresponded to the number of words accurately recalled (Cronbach's α = .83).

**Visuospatial working memory**

*Matrices span -Backward* (MS-B; Cornoldi & Vecchi, 2003). The children had to memorise and recall the positions of black cells that appeared briefly (for 1 second) in different positions on the screen. After a series of black cells had been presented, the children used the mouse to click on the locations where they had seen a black cell appear. The number of cells presented in each series ranged from 2 to 8. The target appeared and disappeared in a visible (4 × 4) grid in the centre of the screen. There was no time limit for recalling the cells in the backward order (Cronbach's α =.79). The score was the number of cells accurately reproduced in the right order.

*Dual tasks visuospatial (DT-VS)* (Mammarella & Cornoldi, 2005). In the DT-VS, a series of two-dimensional 4 × 4 grids, each comprising 16 empty cells, was shown on the screen. Seven of the sixteen cells were always colored in gray while the others were white. The task was administered in sets of three grids, in which a black dot appeared in one of the cells, and then disappeared. The children were asked to press the spacebar if the dot appeared in a gray cell, and also to remember the last position of the dot (in the third grid in each set). The sets of grids were grouped together in gradually longer series, so that the number of dots to remember ranged from 2 to 6 (Cronbach's α = .82). The score corresponded to the number of dot positions accurately recalled in the right order.

*Dot matrix task* (DOT, derived from Miyake, Friedman, Rettinger, Shah, & Hegarty, 2001). In this task, children were shown a matrix equation that they were asked to verify, then a dot appeared in a 5 × 5 grid and they had to remember its position. The matrix equation involved adding or subtracting simple line drawings. After a given series of pairs of equations and grids, the positions of the dots in the various grids had to be recalled by clicking with the mouse on an empty grid. The series increased in length so that from 2 to 5 dot positions had to be remembered. The score corresponded to the number of dot positions correctly recalled (Cronbach's α =.74).

**2.2.2 Mathematical and reading literacy**

*Mathematical literacy* (INVALSI, 2011). The appropriate version of the INVALSI (Italian Institute for the Assessment of the Instruction System) test was used for each school grade. The INVALSI test measures four areas: *space and figures* (MAT-SF) includes problems on two or three-dimensional solids, or other mainly geometry-related tasks; *numbers* (MAT-N) measures number fractions and other mathematical elements; *relations and functions* (MAT-RF) involves solving problems with equivalences or algebraic expressions; and *data and prediction* (MAT-DP) consists of calculating the probability of an event, or means, medians, frequencies, and other statistical properties. The task administered consisted of 31 questions with 49 items for the sixth-graders, and 25 questions with 46 items for the eighth-graders. The proportion of correct answers was considered to make the results comparable. The test took 75 minutes to complete (Cronbach's α = .84 for the sixth-graders and .85 for the eighth-graders).

*Reading literacy* (INVALSI, 2011). The INVALSI included two types of task, covering reading comprehension and grammar. For reading comprehension (READ-RC), the children were shown two or three written texts, and they had to answer several multiple-choice or open-ended questions. For grammar (READ-G), they were asked to answer questions on Italian language spelling, morphology and lexicon. The task consisted of 45 questions for the children in sixth grade (36 on reading comprehension, 9 on grammar), and 49 questions for the eighth-graders (38 on reading comprehension, 11 on grammar). Here again, the proportion of correct answers was considered to make the results comparable. The test lasted 75 minutes (Cronbach's α = .87 for the sixth-graders, and .88 for the eighth-graders).

**2.3 Procedure**

The tasks were chosen by agreement with the schools taking part in the study, and were administered as part of a broader study on the relationship between cognitive, emotional, and academic achievement factors in schoolchildren. Participants were tested in different stages: (a) WM assessment via individual sessions in a quiet room away from the classroom, lasting approximately 30 min.; and (b) academic achievement assessment via two group sessions lasting 75 minutes each. WM tasks were administrated from February to March, while mathematical and reading literacy were administrated on May of the same school year. For the latter, it is worth noting that the instructions in the INVALSI protocol were followed, i.e. an external investigator and a teacher was always present during the tests to minimize cheating or other issues.

During the individual sessions, the WM tasks were also administered in the following order: (1) MS-B; (2) DT-V; (3) DOT; (4) LST; (5) WS-B; (6) DT-VS. All the WM tasks administered during the individual sessions were programmed using E-prime 2 software and presented on a 15-inch touchscreen laptop. After two practice trials, each task began on the easiest level and gradually became more difficult, with two trials for each progressively higher level of complexity. The partial credit scoring method was used for scoring purposes (Conway et al., 2005; Giofrè & Mammarella, 2014).

**3. Results**

**3.1 Data analysis**

The R program (R Core Team, 2014) with the "lavaan" library (Rosseel, 2012) was used. Model fit was assessed using various indexes according to the criteria suggested by Hu and Bentler (1999). We considered the chi-square ($\chi^2$), the comparative fit index (CFI), the non-normed fit index (NNFI), the standardized root mean square residual (SRMR), and the root mean square error of approximation (RMSEA). The chi-square difference ($\Delta\chi^2$), and the Akaike information criterion (AIC) were also used to compare the fit of alternative models (Kline, 2011).

**3.2 Preliminary analyses**

The assumption of multivariate normality and linearity was tested using Mardia's test. The results were not statistically significant, confirming that normality could be assumed. The chosen estimation method (maximum likelihood) is therefore robust against several types of violation of the multivariate normality assumption (Bollen, 1989). A partial correlation analysis was conducted with age in months, taking the children's different ages and school grades into account. Partial correlations were used in all the analyses (Alloway et al., 2006; Giofrè, et al., 2013).

**3.3 WM models**

In the first CFA model, model 1, we hypothesized the presence of a single WM factor, without distinguishing between verbal and visuospatial WM (Figure 1). The model`s fit was not adequate, $\chi2(9) = 50.39$, $p < .001$, $RMSEA = .181$, $SRMR = .076$, $CFI = .846$, $NNFI = .743$, $AIC = 5290$, and we decided to test an alternative model distinguishing between verbal and visuospatial modality.

In model 2, we distinguished between verbal and visuospatial WM (Figure 1). This model had a better fit, $\chi2(8) = 9.92$, $p = .271$, $RMSEA = .041$, $SRMR = .045$, $CFI = .993$, $NNFI = .987$, $AIC = 5251$, and was statically superior compared to the previous one, $\Delta\chi^2(1) = 40.48$, $p < .001$, and was retained for further analyses. In this model, the correlation between WM-V and WM-VS was high (.67, with 95% confidence intervals of .541-.808]) (Cohen, 1988).

**3.4 Overall CFA model**

Model 3. In this model, four factors were considered: mathematics (math), reading (read), WM-V, and WM-VS (Figure 1). All the patterns were statistically significant, and the model had an adequate fit, $\chi2(48) = 67.98$, $p = .030$, *RMSEA* = .054, *SRMR* = .057, *CFI* = .967, *NNFI* = .955, *AIC* = 12361, so the model was retained (Figure 1).

Figure 1 about here

**3.5 Hierarchical regressions**

In the final set of analyses, we used variance partitioning to examine the unique and shared portions of the variance explained in mathematics and reading. This gave us the opportunity to see the specific contribution of WM-V and WM-VS to reading and mathematics. In order to partition the variance, we conducted a series of regression analyses to obtain $R^2$ values from different combinations of predictor variables. For each variable entered into the regression, inter-factor correlations were used (see Giofrè et al., 2014 for a similar procedure). We found that when WM-V and WM-VS were entered simultaneously, a large variance was predicted in both mathematics (40.4%) and reading (31.3%). As shown in Venn Diagrams (Figure 2), we found that that the variance shared by WM-V and WM-VS predicted a conspicuous portion of the variance in both reading and mathematics. In addition, WM-V accounted for a large portion of the unique variance in reading (15%), while WM-VS accounted for a large portion of the unique variance in mathematics (10%), although WM-V also slightly contributed to the variance in mathematics (4%).

Figure 2 about here

**4. Discussion**

The main aim of this paper was to investigate whether verbal and visuospatial WM tasks are distinguishable and whether correlations follow a pattern between a particular domain, with WM-V more strongly related to reading, and WM-VS more strongly related to mathematics.

Concerning our first hypothesis, we found that verbal and visuospatial WM were highly related but distinguishable. The correlation between WM-V and WM-VS was .67, with confidence intervals not including 1, which confirms that WM-V and WM-VS are empirically distinguishable. In fact, about 45% of the variance is shared between these two factors (i.e., domain-general), while about 55% of the variance has domain-specific origin. This latent correlation is in line with other studies of children of various ages and of different countries using a completely different set of tasks (Alloway & Alloway, 2013; Swanson, 2017). This result is also in accordance with neuroimaging studies indicating that different brain regions are activated during verbal and visuospatial tasks (see Donolato, Giofrè, & Mammarella, 2017 for a review) and by evidence indicating different developmental trajectories of verbal and visuospatial WM (Alloway & Alloway, 2013; Swanson, 2017).

According to our second hypothesis, we found that verbal and visuospatial WM had a different predictive value on mathematics and reading. In agreement with previous findings, taken together verbal and visuospatial WM predicted a large portion of the variance in both mathematics (40.4%) and reading (31.3%) (see for example Swanson & Alloway, 2012 for a similar finding). This result is in line with a large body of research strengthening the importance of cognitive control (or executive control) over reading and mathematics (e.g., Engle, Tuholski, Laughlin, & Conway, 1999). However, we also found that verbal and spatial WM showed some domain-dependence on other higher cognitive skills (Friedman & Miyake, 2000; Handley, Capon, Copp, & Harper, 2002; Shah & Miyake, 1996). In fact, the variance explained in mathematics and reading by verbal and visuospatial WM is not trivial. Regarding mathematics, we found that visuospatial WM accounted for a large portion of the unique variance in mathematics (10%), with a smaller portion of the variance uniquely explained by verbal WM (about the 4%). This finding is also consistent with previous research indicating the importance of visuospatial WM in mathematics (Caviola et al., 2014; Mammarella et al., 2017; Szűcs, 2016; Szűcs et al., 2013). As for reading, we found that an important portion of the variance is uniquely explained by verbal WV (about 15%), which is in agreement with previous research strengthening the importance of verbal WM in reading (Carretti, Borella, Cornoldi, & De Beni, 2009; Oakhill et al., 2005; Vellutino, 2004).

The present study has some limitations that should be addressed by future research. First, it considered only a limited sample of typically-developing 11- and 13-year-old children, so our findings can only be applied to this particular school population. Moreover, the sample size considered in the present study prevent us from further investigating possible age related effects, which should therefore considered in future studies (preferably longitudinal). Second, we found that mathematics is more influenced by visuospatial WM and reading is more influenced by verbal WM. This is not entirely consistent with other research findings that indicate a greater importance of domain-general resources over domain-specific resources (Swanson, 2017). Importantly, we considered verbal and spatial WM, reading and mathematics at the latent level. Moreover, reading vocabulary and mathematic computations were used as measures of achievement in Swanson (2017), which strongly differed from our mathematics and reading achievement measures used in the present study. However, future research on this topic is needed to shed further light on domain-general and domain specific individual contributions to mathematics and reading.

The present findings have interesting clinical and educational implications. For example, in children with specific learning disabilities, WM is consistently impaired apparently suggesting that children with specific reading or mathematics learning disorders could potentially benefit from either verbal or visuospatial WM trainings (Dunning, Holmes, & Gathercole, 2013). A recent meta-analysis supported the efficacy of specific WM trainings to determine several improvements in a wide range of untrained WM tasks. However, these gains do not clearly translate into enhancements in academic progress, including reading or mathematics (Melby-Lervåg, Redick, & Hulme, 2016). In terms of educational implications, typically developing children would benefit from a reduced load on their verbal or visuospatial WM (Gathercole et al., 2006), and also benefit from consolidating the acquisition of simpler reading or mathematics skills before moving on to more complex ones.

To sum up, in agreement with a large body of literature, we found verbal and visuospatial WM tasks can be distinguished. In fact, the correlation between verbal and visuospatial WM was high but not substantial, indicating that more than half of the variance is not shared between the two modalities. As for the predictive utility of verbal and visuospatial WM, we found that the unique contribution of verbal WM is larger in reading and smaller in mathematics, and vice versa for visuospatial WM, in which the unique contribution is substantial in mathematics and very small in reading.

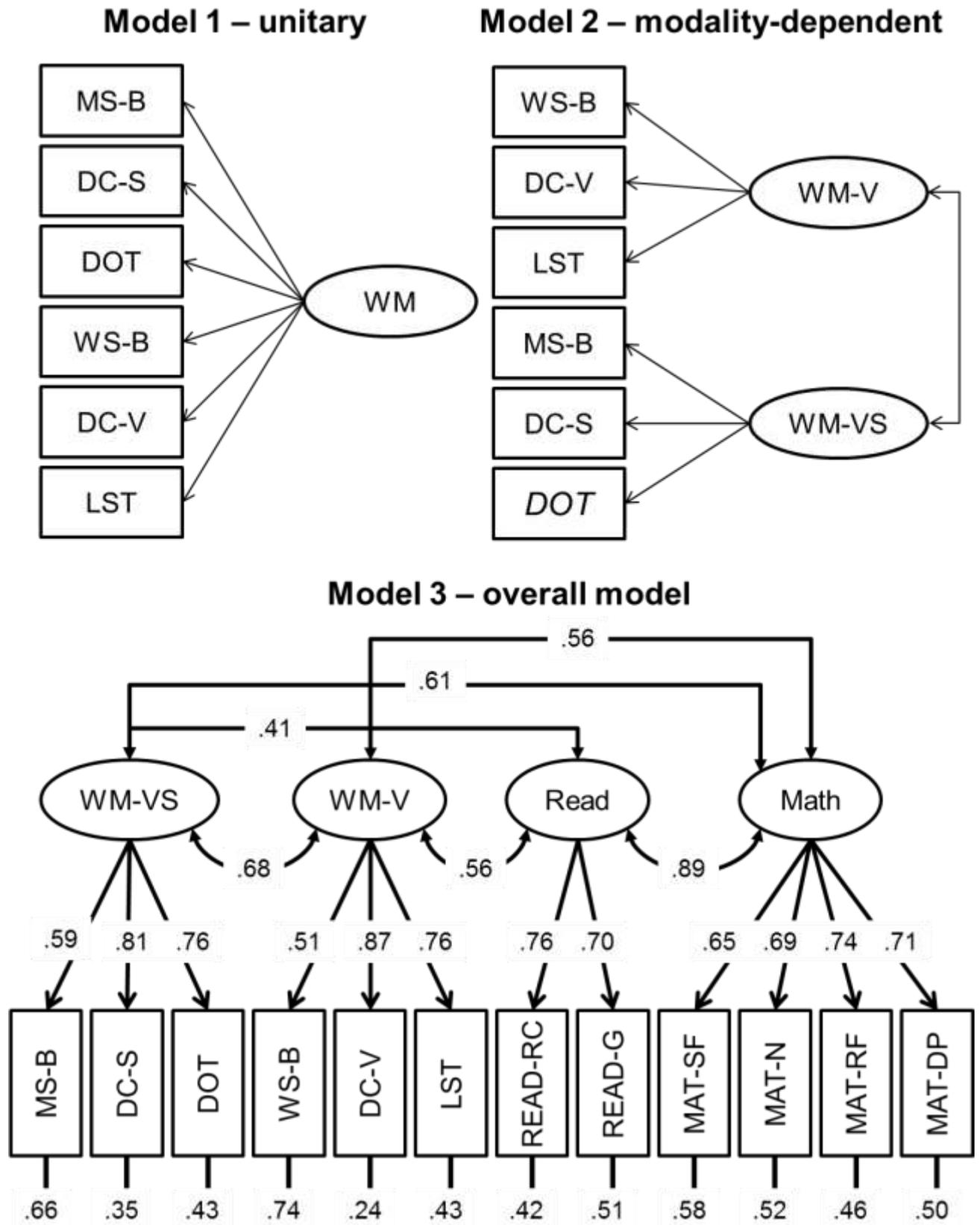

*Figure 1*. Conceptual diagrams, and Measurement CFA model. All paths are significant at .05 level. MS-B = Matrix Span Backward, DT-S = Dual Task Spatial, DOT = Dot-matrix, WS-B = Word Span Backward, DT-V = Dual Task Verbal, LST = Listening span test, -SF = Space and Figures, -N = Numbers, -RF = relations and functions, DP = data and previsions, RC = reading comprehension, G = Grammar, Math = Mathematics literacy, Reading = Reading literacy.

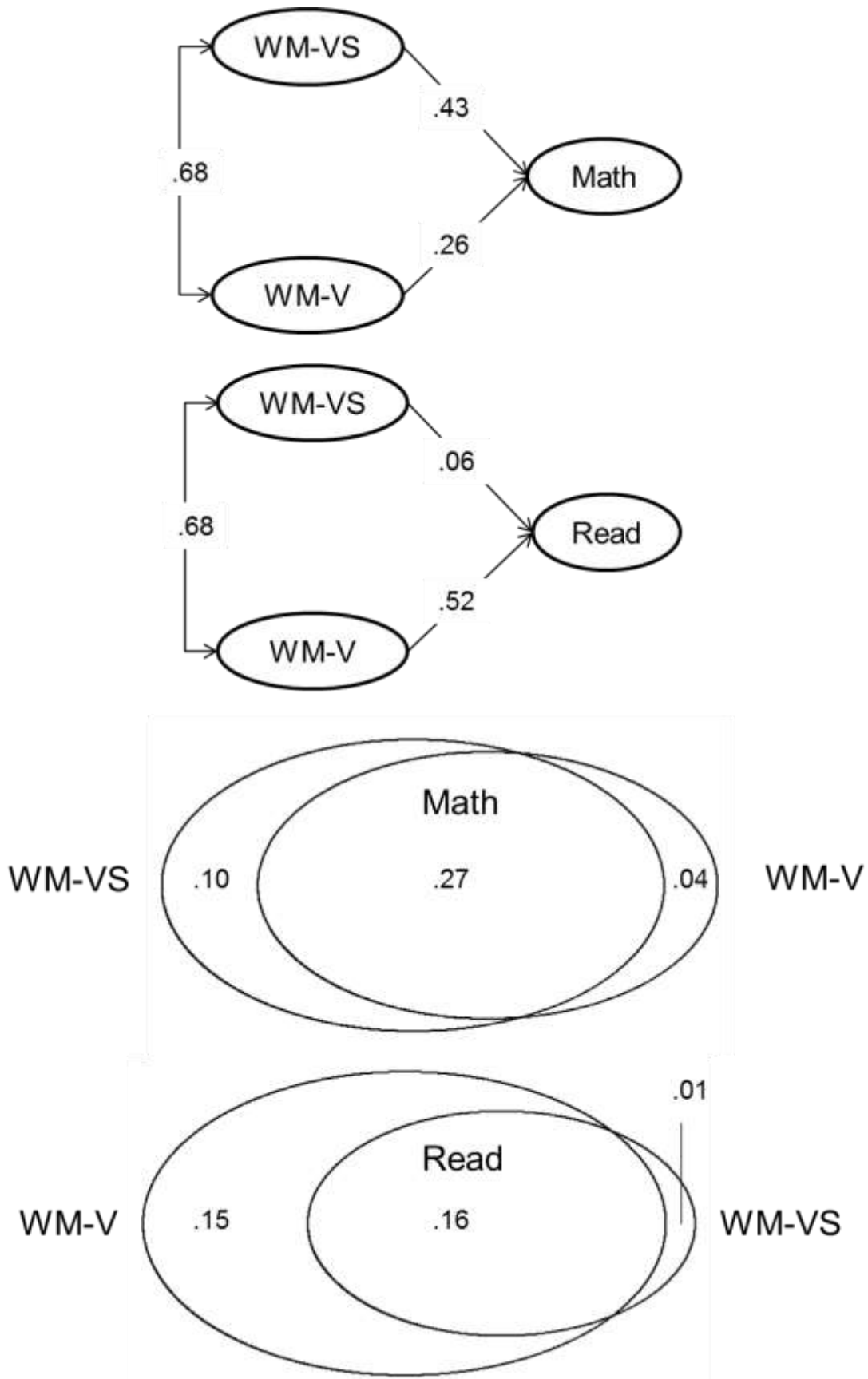

*Figure 2.* **Top:** Regressions with math and reading as response variables, and WM-V and WM-VS as predictors. **Bottom:** Venn diagrams indicating the shared and unique variance explained in mathematics and reading by verbal (WM-V) and visuospatial (WM-VS) working memory.